\documentclass[12pt]{article}%
\usepackage{amssymb}
\usepackage{amsfonts}
\usepackage{amsmath}
\usepackage{graphicx}%
\providecommand{\U}[1]{\protect\rule{.1in}{.1in}}
\newtheorem{theorem}{Theorem}

\newtheorem{lem}{Lemma}
\newtheorem{lemma}[lem]{Lemma}
\newtheorem{rem}{Remark}
\newtheorem{remark}[rem]{Remark}

\begin{document}
\begin{align*}
&  \text{{\LARGE \ Comments on the height reducing property }}\\
&  \text{{\LARGE \ }}%
\end{align*}

\begin{center}
Shigeki AKIYAMA \ \ and \ \ \bigskip Toufik ZAIMI
\end{center}

ABSTRACT.\textbf{\ \ }\textit{A complex number }$\alpha$\textit{\ is said to
satisfy the height reducing property if there is a finite subset, say }$F,$
\textit{of the ring }$\mathbb{Z}$ \textit{of the rational integers\ such that}
$\ \mathbb{Z}[\alpha]=F[\alpha].$\textit{\ This problem has been considered by
several authors, especially in contexts related to self affine tilings, and
expansions of real numbers in non-integer bases. We continue, in this paper,
the description of the numbers satisfying the height reducing property, and we
specify a related characterization of the roots of integer polynomials with
dominant term. }

\bigskip\medskip

\begin{center}
\textbf{1. Introduction}
\end{center}

For a subset $F$ of the complex field $\mathbb{C,}$ and for $\alpha
\in\mathbb{C,}$ we denote by $F[\alpha]$ the set of polynomials with
coefficients in $F,$ evaluated at $\alpha,$ i. e.,
\[
F[\alpha]=\{%
{\displaystyle\sum\limits_{j=0}^{n}}
\varepsilon_{j}\alpha^{j}\mid(\varepsilon_{0},...,\varepsilon_{n})\in
F^{n+1},\text{ }n\in\mathbb{N}\},
\]
where $\mathbb{N}$ is the set of non-negative rational integers. In
particular, when $F$ \ is the ring $\mathbb{Z}$ of the rational integers, the
set $F[\alpha]$ is the $\mathbb{Z}-$module generated by the integral powers of
$\alpha.$ It is well known that there is $N\in\mathbb{N}$ such that
$\mathbb{Z}[\alpha]=\{\varepsilon_{0}+\cdot\cdot\cdot+\varepsilon_{N}%
\alpha^{N}\mid(\varepsilon_{0},...,\varepsilon_{N})\in\mathbb{Z}^{N+1}\}$ if,
and only if, $\alpha$ is an algebraic integer; moreover, the smallest possible
value for $N,$ in this case, is $\deg(\alpha)-1,$ where $\deg(\alpha)$ is the
degree of $\alpha$ \cite{Nar}.

The analog height reducing problem, for the ring $\mathbb{Z}[\alpha],$ which
consists in the existence of a set, say again $F,$ satisfying
\begin{equation}
F\subset\mathbb{Z},\text{ }\mathbb{Z}[\alpha]=F[\alpha]\text{ \ and }F\text{
\ finite,}%
\end{equation}
----------------------------------------------------------------------------------------------



\newpage has been considered by several authors (see the references in
\cite{ADJ,Dub}), especially in contexts related to self affine tilings, and
expansions of real numbers in non-integer bases. A result of Lagarias and
Wang, cited in \cite{ADJ,Dub}, implies that an expanding algebraic integer
$\alpha,$ that is an algebraic integer whose conjugates are of modulus greater
than one, satisfies (1) with $F=\{0,\pm1,...,\pm(\left\vert Norm(\alpha
)\right\vert -1)\}.$ Recently, Akiyama, Drungilas and Jankauskas obtained a
direct proof of this last mentioned result, but with a greater finite set $F$
\ \cite{ADJ}. It is worth noting that Proposition 3.1 of \cite{FroSte} yields
to the same conclusion. Also, Lemma 1 of \cite{ADJ} asserts that an algebraic
integer, with modulus greater than 1, satisfying the height reducing property,
is an expanding algebraic integer. Next we continue the description of numbers
satisfying (1).

\bigskip

\begin{theorem}
\label{HRP} Let $\alpha\in\mathbb{C}$. Then, the following propositions are true.

\begin{enumerate}
\item If $\alpha$ satisfies the height reducing property , then $\alpha$ is an
algebraic number whose conjugates are all of modulus $1,$ or all of modulus
greater than $1$.

\item If $\alpha$ is a root of unity, or an algebraic number whose conjugates
are of modulus greater than $1,$ then $\alpha$ satisfies the height reducing property.
\end{enumerate}
\end{theorem}

It is clear, by Kronecker's theorem (see for instance \cite{Nar}), that an
algebraic integer whose conjugates belong to the unit circle is a root of
unity. To obtain a characterization of the numbers which satisfy (1), it
remains to consider the case where the conjugates of the algebraic number
$\alpha$ belong to the unit circle, and are not roots of unity. In this last
situation the minimal polynomial $M_{\alpha}$ of $\alpha$ is reciprocal, i.
e., $M_{\alpha}(x)=x^{\deg(M_{\alpha})}M_{\alpha}(1/x),$ $\deg(M_{\alpha})$
(which is equal to $\deg(\alpha))$ is even, and the greatest number, say
$m(\alpha),$ of conjugates of $\alpha$ which are multiplicatively independent
(see the definition in Lemma \ref{KroA} below) satisfies the relation $1\leq
m(\alpha)\leq\deg(\alpha)/2,$ since the roots of $M_{\alpha}$ are pairwise
complex conjugates and $\arg(\alpha)/\pi\notin\mathbb{Q}$ (i.e., $\alpha$ is
not a root of unity).

\bigskip

\begin{theorem}
\label{Circle} Let $\alpha$ be an algebraic number whose all conjugates lie on
the unit circle. If $m(\alpha)\ge\deg(\alpha)/2-1,$ or $m(\alpha)=1$, then
$\alpha$ satisfies the height reducing property.
\end{theorem}

\begin{remark}
\label{LowDeg} It follows immediately from Theorem \ref{Circle} that $\alpha$
satisfies the height reducing property when $\deg(\alpha)\leq6$. We expect
that height reducing property holds for any algebraic $\alpha$ whose
conjugates lie on the unit circle. However we find two examples of degree 12
that none of our methods apply in the Appendix.
\end{remark}

\begin{remark}
\label{comp} There is an algorithm to determine $m(\alpha)$. In fact if
$\alpha_{1},\dots\alpha_{m}$ are multiplicatively dependent, then Lemma 4.1 in
Waldschmidt \cite{Wal} gives an explicit upper bound $B$ so that the equation
$\prod_{i=1}^{m}\alpha_{i}^{k_{i}}=1$ has a non-trivial solution
$(k_{1},...,k_{m})\in(\mathbb{Z}\cap\lbrack-B,B])^{m}$. However the bound $B$
is too large to examine. We employ Lemma 3.7 of de Weger \cite{deW} to reduce
this bound by LLL algorithm. Details and numerical results will be shown in
the Appendix.
\end{remark}

Following \cite{Dub}, we say that a non-zero polynomial $P=P(x)=c_{0}%
+\cdot\cdot\cdot+c_{\deg(P)}x^{\deg(P)}\in\mathbb{C}[x]$ has a dominant term
(resp., has a dominant constant term) if there is $k\in\{0,...,\deg(P)\}$ such
that $\left\vert c_{k}\right\vert \geq%
{\displaystyle\sum\limits_{j\neq k}}
\left\vert c_{j}\right\vert $ (resp., such that $\left\vert c_{0}\right\vert
\geq%
{\displaystyle\sum\limits_{1\leq j}}
\left\vert c_{j}\right\vert ).$ In connection with a property studied by
Frougny and Steiner \cite{FroSte}, about minimal weight expansions, Dubickas
obtained recently \cite{Dub}, some characterizations of complex numbers which
are roots of integer polynomials (i. e., polynomials with rational integer
coefficients) having a dominant term.

\bigskip

\textbf{Theorem A (\cite{Dub})} \textit{Let }$\alpha\in\mathbb{C}%
.$\textit{\ Then, the following assertions are true.}

\textit{(i) The number }$\alpha$\textit{\ is a root of \ an integer polynomial
with dominant term if, and only if, }$\alpha$\textit{\ is a root of unity, or
}$\alpha$\textit{\ is an algebraic number without conjugates of modulus }$1.$

\textit{(ii) The number }$\alpha$\textit{\ is a root of an integer polynomial
with dominant constant term if, and only if, }$\alpha$\textit{\ is a root of
unity, or }$\alpha$\textit{\ is an algebraic number all of whose conjugates
are of modulus greater than }$1.$

\bigskip

The other aim of this paper is to show two simple generalizations of Theorem
A. The first one is an integral version of Theorem A (i). To state the second
one, let us introduce the following "definition-precision" : \textit{We say
that the non-zero polynomial }$P,$\textit{\ defined above,\ has a }$k-
$\textit{th} \textit{dominant term, }(\textit{resp., has a }$k-$\textit{th
strictly dominant term}),\textit{\ where }$k\in\{0,...,\deg(P)\},$\textit{\ if
}$\left\vert c_{k}\right\vert \geq%
{\displaystyle\sum\limits_{j\neq k}}
\left\vert c_{j}\right\vert $\textit{\ }(\textit{resp., if \ }$\left\vert
c_{k}\right\vert >%
{\displaystyle\sum\limits_{j\neq k}}
\left\vert c_{j}\right\vert ).$ \textit{The polynomial }$P$\textit{\ has a
strictly dominant term, when it } \textit{has some }$k-$\textit{th}
\textit{strictly dominant term.}

\bigskip

\begin{theorem}
\label{Dominant} Let $\alpha\in\mathbb{C}.$ Then, the following propositions
are true.

\begin{enumerate}
\item The number $\alpha$ is a root of an (resp., of a monic) integer
polynomial with $k$-th dominant term if, and only if, $\alpha$ is a root of
unity, or $\alpha$ is an algebraic number (resp., algebraic integer)\ having
at most $k$ conjugates inside the unit disk and no conjugates on the unit circle.

\item The number $\alpha$ is a root of an (resp., of a monic) integer
polynomial with $k$-th strictly dominant term if, and only if, $\alpha$ is an
algebraic number (resp., algebraic integer) having at most $k$ conjugates
inside the unit disk and no conjugates on the unit circle.
\end{enumerate}
\end{theorem}

Obviously, Theorem A (ii) is a corollary of Theorem \ref{Dominant} (i), with
$k=0.$ Theorem \ref{Dominant} (i) implies Theorem A (i), too. It follows also
from Theorem 3 (ii) that a complex number is a root of some (resp., some
monic) integer polynomial\ with strictly dominant term\ if, and only if,
it\ is an algebraic number (resp., algebraic integer) without\ conjugates on
the unit circle.

In these pages when we speak about conjugates, norm, minimal polynomial and
degree of an algebraic number we mean over the field of the rationals
$\mathbb{Q}.$ A unit is an algebraic integer whose norm is $\pm1.$The proofs
of the theorems above appear in the last section. Theorem \ref{Dominant} of
the present manuscript, and some parts of the proofs of Lemmas 1 and 6 of
\cite{ADJ} are used to show Theorem \ref{HRP}. Lemmas 5 and 6 of \cite{Dub}
are the main tool of the proof of Theorem \ref{Dominant}; these lemmas,
together with some auxiliary results we need to prove Theorem \ref{Circle},
are exhibited in the next section.

\bigskip

\begin{center}
\bigskip\textbf{\bigskip2. Some lemmas}
\end{center}

The following result is the main tool of the first part of the proof of
Theorem \ref{Circle}.

\smallskip

\begin{lemma}
\label{KroA} Let $\alpha_{1},\dots,\alpha_{m}$ be conjugates, with modulus
one, of an algebraic number $\alpha.$ Assume that $\alpha_{1},\dots,\alpha
_{m}$ are multiplicatively independent, i.e., any equation of the form
$\prod_{j=1}^{m}\alpha_{j}^{k_{j}} =1$ where $(k_{1},\dots,k_{m})\in
\mathbb{Z}^{m},$ implies $(k_{1},\dots,k_{m})=(0,\dots,0).$ Then for any
$\varepsilon>0,$ there is a positive rational integer $K=K(\alpha
,m,\varepsilon)$ such that for any non-zero complex numbers $\beta_{1}%
,\dots,\beta_{m}$ there is a non-negative rational integer $l\leq
K$\textit{\ \ satisfying }$|\arg(\beta_{j}\alpha_{j}^{l})|\leq\varepsilon
,$\textit{\ }$\forall$\textit{\ }$j\in\{1,\dots,m\}.$
\end{lemma}

\textbf{Proof.} The existence of the constant $K,$ satisfying the above
mentioned condition, is a corollary of a quantitative version of Kronecker's
approximation theorem due to Mahler \cite{Mah} (c.f. Vorselen \cite{Vor}). The
necessary assumption of the lower bound follows from Baker's theory of linear
forms in logarithms (see \cite{Bak,BakWus}).\hfill$\Box$

\smallskip

To simplify the computation in the proof of Theorem \ref{Circle}, let us show
the following lemma.

\begin{lemma}
\label{L2} Let $z$\textit{\ and }$w$ be complex numbers satisfying $z\neq
0,$\textit{\ }$\left\vert \arg(z)\right\vert \leq2\pi/5$\textit{\ and
}$\left\vert w\right\vert \leq1.$ Then for any real number $r$ $\in
(0,4\left\vert z\right\vert /145),$ we have
\[
|z+r(w-5)|<\left\vert z\right\vert .
\]

\end{lemma}

\textbf{Proof. }Set $z:=\delta\exp(i\theta),$ $w:=\rho\exp(i\phi)$ and
$(z+r(w-5))\exp(-i\theta):=a+ib,$ where $i^{2}=-1,$ $\{\delta,$ $\theta,$
$\rho,$ $\phi,$ $a,$ $b\}\subset\mathbb{R}$ and $\mathbb{R}$ is the real
field. Then $a=\delta+r\rho\cos(\phi-\theta)-5r\cos(\theta),$ $b=r\rho
\sin(\phi-\theta)+5r\sin(\theta),$ $0<\delta-6r\leq a\leq\delta-(5\cos
(2\pi/5)-1)r\leq\delta-r/2,$ \ $\ |b|\leq6r$ and so
\[
|z+r(w-5)|\leq\sqrt{(\delta-r/2)^{2}+36r^{2}}<\delta\text{.\hfill$\Box$}
\]

\bigskip

\begin{lemma}
\label{Int} Let $\alpha$ be an algebraic number of degree $d$. Then
$\mathbb{Z}[\alpha] \cap\mathbb{Z}[1/\alpha]$ is an \textit{order}, i.e., the
subring of the integer ring of $\mathbb{Q}(\alpha)$ sharing the identity, as
well as a free $\mathbb{Z}$-submodule of rank $d$.
\end{lemma}

\textbf{Proof.} Put $\mathcal{O}=\mathbb{Z}[\alpha]\cap\mathbb{Z}[1/\alpha]$.
If $\alpha$ is an algebraic integer, then we have $\mathbb{Z}[\alpha
]\subset\mathbb{Z}[1/\alpha]$ and the statement is trivial. Assume that
$\alpha$ is not an algebraic integer, and take an ideal $\mathfrak{p}$ which
divides the denominator of the fractional ideal $(\alpha)$. Then the
denominator of the principal ideal $(x)$ for $x\in\mathcal{O}$ is not
divisible by $\mathfrak{p}$. This shows that every element of $\mathcal{O}$ is
an algebraic integer and $\mathcal{O}$ is a $\mathbb{Z}$-module of rank not
greater than $d$. Denote by $\sum_{n=0}^{d}c_{n}x^{n}$ the minimal polynomial
of $\alpha$. Then from the relation
\[
c_{d}\alpha=-\sum_{n=0}^{d-1}c_{n}\alpha^{n-d+1}\in\mathbb{Z}[1/\alpha],
\]
and the fact that $c_{d}\alpha$ is an algebraic integer, we see that
\[
\mathbb{Z}[c_{d}\alpha]\subset\mathbb{Z}[\alpha]\cap\mathbb{Z}[1/\alpha].
\]
This shows that the rank of $\mathcal{O}$ is not less than $d$. \hfill$\Box$
\bigskip

\begin{lemma}
\label{L6} Let $\alpha$ be an algebraic number of degree $2d$ whose all
conjugates are of modulus one. Let $\alpha_{j}\ (j=1,\dots,d)$ be the
conjugates of $\alpha$ lying in the upper half plane. If $m(\alpha)=d-1$ then
there is a vector $(a_{1},\dots a_{d})\in\{-1,1\}^{d}$ and a root of unity
$\zeta$ such that $\prod_{j=1}^{d}\alpha_{j}^{a_{j}}=\zeta$.
\end{lemma}

\textbf{Proof}\textit{.} If $m(\alpha)=0$ then $\alpha=\pm i$ and $\alpha$ is
a root of unity. Suppose $m(\alpha)\geq1.$ Then $d\geq2,$ and by
$m(\alpha)=d-1$, there is  $(b_{1},\dots,b_{d})\in\mathbb{Z}^{d}%
\setminus\{(0,\dots,0)\}$ such that $\prod_{j=1}^{d}\alpha_{j}^{b_{j}}=1$. It
suffices to show that there exists a positive rational integer $b$ satisfying
$|b_{j}|=b$ for all $j$. If not, then we may assume that $|b_{1}|>|b_{2}%
|=\min_{j=1}^{d}|b_{j}|$. Applying a conjugate map $\sigma$ which sends
$\alpha_{2}$ to $\alpha_{1}$, \ we obtain $\prod_{j=1}^{d}\alpha_{j}^{c_{j}%
}=1,$ with $c_{1}=b_{2}$, and so
\[
\prod_{i=2}^{d}\alpha_{i}^{b_{1}c_{i}-b_{2}b_{i}}=1.
\]
Since $|c_{j}|=|b_{1}|$ for some $j$, this last multiplicative relation is non
trivial, and yields, together with the equation $\prod_{j=1}^{d}\alpha
_{j}^{b_{j}}=1,$ to the inequality $m(\alpha)<d-1$. \hfill$\Box$ \bigskip

\begin{lemma}
\label{MinusOne} Let $\alpha$ be an algebraic number of degree $2d\geq6$ whose
all conjugates are of modulus one. Let $\alpha_{j}\ (j=1,\dots,d)$ be the
conjugates of $\alpha$ lying in the upper half plane. If $m(\alpha)=d-1$ then
there is a positive integer $K=K(\alpha)$ such that for any non-zero complex
numbers $\beta_{1},\dots,\beta_{d}$ there is a non-negative integer $\ell\leq
K$ such that $|\arg(\beta_{j}\alpha_{j}^{\ell})|\leq2\pi/5$ for $j=1,\dots,d$.
\end{lemma}

\textbf{Proof.}  Lemma \ref{L6} asserts that there is a positive rational
integer $b$ such that
\[
\alpha_{1}^{b}=\alpha_{2}^{\pm b}\dots\alpha_{d}^{\pm b}%
\]
for a fixed choice of $\pm$'s and $\alpha_{2}^{b},\dots,\alpha_{d}^{b}$ are
multiplicatively independent. So substituting $\alpha_{j}^{\pm b}$ to $\alpha_{j}$ 
for each $j,$ we may assume that
\[
\alpha_{1}=\alpha_{2}\dots\alpha_{d}.
\]
This implies
\begin{equation}
\beta_{1}\alpha_{1}^{\ell}=\beta_{1}(\prod_{j=2}^{d}\beta_{j}\alpha_{j}^{\ell
})/(\prod_{j=2}^{d}\beta_{j})\label{Dep}%
\end{equation}
for any $\ell$. Fix a small $0<\varepsilon<\pi/15$ and apply Kronecker's
approximation theorem as in Lemma \ref{KroA} to the following three sets
of\ $(d-1)$ inequalities:

\begin{itemize}
\item $\left|  \arg(\beta_{2} \alpha_{2}^{\ell})- \frac{\pi}3\right|
<\varepsilon,\ \left|  \arg(\beta_{3} \alpha_{3}^{\ell})- \frac{\pi}3\right|
<\varepsilon,\ \left|  \arg(\beta_{j} \alpha_{j}^{\ell})\right|
<\varepsilon\ (j\ge4) $

\item $\left|  \arg(\beta_{2} \alpha_{2}^{\ell})+ \frac{\pi}3\right|
<\varepsilon,\ \left|  \arg(\beta_{3} \alpha_{3}^{\ell})- \frac{\pi}3\right|
<\varepsilon,\ \left|  \arg(\beta_{j} \alpha_{j}^{\ell})\right|
<\varepsilon\ (j\ge4) $

\item $\left|  \arg(\beta_{2} \alpha_{2}^{\ell})+ \frac{\pi}3\right|
<\varepsilon,\ \left|  \arg(\beta_{3} \alpha_{3}^{\ell})+ \frac{\pi}3\right|
<\varepsilon,\ \left|  \arg(\beta_{j} \alpha_{j}^{\ell})\right|
<\varepsilon\ (j\ge4) $
\end{itemize}

Then we can find a common $K=K(\alpha)$ such that these 3 systems are solvable
with $\ell_{j}\leq K\ (j=1,2,3)$. We see from (\ref{Dep}) that one of the
systems gives the solution of our problem. \hfill$\Box$ \bigskip

\begin{lemma}
\label{Dep2} Let $\alpha_{1},\alpha_{2}$ be two conjugates of an algebraic
number $\alpha$. Assume that $\alpha$ is not a unit and there is
$(a,b)\in\mathbb{Z}^{2}\setminus\{(0,0)\}$ with $\alpha_{1}^{a}\alpha_{2}%
^{b}=1$. Then $|a|=|b|$.
\end{lemma}

\textbf{Proof.} By the prime ideal decomposition of the fractional ideals
$(\alpha_{1})$ and $(\alpha_{2})$ in the minimum decomposition field of
$\alpha$, we have
\[
(\prod_{j=1}^{s}\mathfrak{p}_{j}^{ae_{j}})(\prod_{j=1}^{s}\mathfrak{p}%
_{j}^{be_{j}^{\prime}})=(1)
\]
and so $ae_{j}+be_{j}^{\prime}=0$ for each $j$. If $|a|<|b|$, then
$|e_{j}|>|e_{j}^{\prime}|$ for all $j$, and we claim that this is impossible.
Indeed, consider an index $l$ with $|e_{l}|=\max_{1\leq j\leq s}|e_{j}|$. As
there is a conjugate map which sends $(\alpha_{1})$ to $(\alpha_{2})$, there
exists an index $k$ such that $e_{k}^{\prime}=e_{l},$ and the inequality
$|e_{k}|>|e_{k}^{\prime}|$ leads immediately to a contradiction. \hfill$\Box$
\bigskip

The following result is the first proposition of Lemma 5 of \cite{Dub}.

\begin{lemma}
\label{Order}\textbf{ (\cite{Dub})} Let $P\in\mathbb{R}[x]$ with dominant
term, and let $\alpha$ be a root of $P$ having modulus one. Then $\alpha$ is a
root of unity.
\end{lemma}

From the proof of Lemma 6 of \cite{Dub}, we easily deduce the following assertion.

\begin{lemma}
\label{L4} (\cite{Dub})\textbf{ } Let $\alpha$ be an algebraic number, having
$l\geq0$ conjugates with modulus less than $1$ and no conjugates on the unit
circle. Then, there is $N\in\mathbb{N}$ such that the polynomial $%
{\displaystyle\prod\limits_{1\leq j\leq\deg(\alpha)}}
(x-\alpha_{j}^{N}),$\textit{\ where }$\alpha_{1},...,\alpha_{\deg(\alpha)}$
are the conjugates of $\alpha,$ has an $l-$th strictly dominant term.
\end{lemma}

\begin{center}
\bigskip

\textbf{3. The proofs}
\end{center}

\textbf{Proof of Theorem \ref{HRP}.} (i) With the notation above, assume that
$\alpha$ satisfies (1) with some finite set $F.$ Then, $F\neq\varnothing$ and
the relation $N\in F[\alpha],$ where $N$ $\in\mathbb{N\cap(}m,\infty) $ and
$m:=\max\{\left\vert \varepsilon\right\vert ,$ $\varepsilon\in F\},$ gives
immediately that $\alpha$ is an algebraic number. Let $\beta$ be a conjugate
of $\alpha.$ Then $\left\vert \beta\right\vert \geq1,$ since otherwise any
element of the set $\mathbb{N\cap(}\frac{m}{1-\left\vert \beta\right\vert
},\infty)$ does not belong to $F[\alpha].$ Now, suppose that $\left\vert
\beta\right\vert =1,$ we have to show that the conjugates of $\alpha$ lie on
the unit circle. If $\deg(\alpha)=1,$ then $\alpha=\pm1$ and the result is
true. Assume that $\deg(\alpha)\geq2.$ Then, the complex conjugate
$\overline{\beta}$ of $\beta$ is also a conjugate of $\alpha,$ and so the
minimal polynomial $M_{\alpha}$ of $\alpha$ divides in the ring $\mathbb{Z}%
[x]$ the polynomial $M_{\alpha}^{\ast}(x):=x^{\deg(\alpha)}M_{\alpha}(1/x),$
as $\overline{\beta}=1/\beta.$ Moreover, the equation $\deg(M_{\alpha}^{\ast
})=\deg(M_{\alpha}),$ yields $M_{\alpha}^{\ast}(x)=cM_{\alpha}(x)$ for some
$c\in\mathbb{Z}$ (in fact we have $c=1$) and so $1/\gamma$ is a conjugate of
$\alpha$ when $\gamma$ is so; thus $\left\vert \gamma\right\vert =1$ since
otherwise one of the numbers $\gamma$ and $1/\gamma$ has modulus less than 1,
and by the above this leads to a contradiction.

(ii) It is clear when $\alpha$ is an $N-$th root of unity, where
$N\in\mathbb{N}^{\ast}:=\mathbb{N}\cap\lbrack1\mathbb{,\infty),}$ that any sum
of the form $%
{\displaystyle\sum\limits_{j=0}^{s}}
a_{j}\alpha^{j},$ where $a_{j}\in\mathbb{Z}$ and $s\in\mathbb{N},$ may be
written
\[%
{\displaystyle\sum\limits_{j=0}^{s}}
\varepsilon_{j}(%
{\displaystyle\sum\limits_{k=1}^{\left\vert a_{j}\right\vert }}
\alpha^{kN})\alpha^{(j+%
{\displaystyle\sum\limits_{l=0}^{j-1}}
\left\vert a_{l}\right\vert N)},
\]
where $\varepsilon_{j}=\mathrm{sgn}(a_{j}),$ and so $\{0,\pm1\}[\alpha
]=\mathbb{Z}[\alpha].$

Now suppose that $\alpha$ is an algebraic number whose conjugates are of
modulus greater than 1. Then Theorem \ref{Dominant} (ii) shows that $\alpha$
is a root of some polynomial $C(x)=c_{0}+c_{1}x+\cdot\cdot\cdot+c_{d}x^{d}%
\in\mathbb{Z}[x],$ with $c_{d}\neq0$ and
\[
\left\vert c_{0}\right\vert >%
{\displaystyle\sum\limits_{j=1}^{d}}
\left\vert c_{j}\right\vert .
\]
Let $R\in\mathbb{Z}[x].$ To prove the relation $R(\alpha)\in F[\alpha],$
where
\[
F:=\{0,\pm1,...,\pm(\left\vert c_{0}\right\vert -1)\},
\]
suppose first that $\deg(R)\in\{0,...,d-1\}.$ Then, $R(x)=A_{0}+\cdot
\cdot\cdot+A_{d-1}x^{d-1},$ for some $(A_{0},...,$ $A_{d-1})\in\mathbb{Z}%
^{d},$ and similarly as in the proof of Theorem 4 of \cite{ADJ}, it suffices
to show, when $A_{0}\notin F,$ \ that
\begin{equation}
R(\alpha)=\varepsilon+\alpha(a_{0}+\cdot\cdot\cdot+a_{d-1}\alpha^{d-1}),
\end{equation}
where $\varepsilon\in F,$ $(a_{0},...,a_{d-1})\in\mathbb{Z}^{d}$ and $%
{\displaystyle\sum\limits_{j=0}^{d-1}}
\left\vert a_{j}\right\vert <%
{\displaystyle\sum\limits_{j=0}^{d-1}}
\left\vert A_{j}\right\vert .$ Since $\left\vert A_{0}\right\vert
\geq\left\vert c_{0}\right\vert ,$ we see that $\left\vert A_{0}\right\vert
=q\left\vert c_{0}\right\vert +\varepsilon,$ for some $q\in\mathbb{N}^{\ast}$
and $\varepsilon\in\mathbb{N}\cap F.$ It follows by the equation $c_{0}%
=-c_{1}\alpha-\cdot\cdot\cdot-c_{d}\alpha^{d},$ that
\[
A_{0}\mathrm{sgn}(A_{0})=qc_{0}\mathrm{sgn}(c_{0})+\varepsilon=\varepsilon
-(qc_{1}\alpha+\cdot\cdot\cdot+qc_{d}\alpha^{d})\mathrm{sgn}(c_{0}),
\]
and so%

\[
A_{0}+\cdot\cdot\cdot+A_{d-1}\alpha^{d-1}=\mathrm{sgn}(A_{0})\varepsilon
+\alpha(a_{0}+\cdot\cdot\cdot+a_{d-1}\alpha^{d-1}),
\]
where $a_{d-1}=-\mathrm{sgn}(c_{0})\mathrm{sgn}(A_{0})qc_{d}$ \ and
$a_{j}=A_{j+1}-\mathrm{sgn}(c_{0})\mathrm{sgn}(A_{0})qc_{j+1}$ for all
$j\in\{0,...,d-2\}.$ Moreover, we have $\mathrm{sgn}(A_{0})\varepsilon\in
F=-F,$ and
\[%
{\displaystyle\sum\limits_{j=0}^{d-1}}
\left\vert a_{j}\right\vert \leq q(%
{\displaystyle\sum\limits_{j=1}^{d}}
\left\vert c_{j}\right\vert )+%
{\displaystyle\sum\limits_{j=1}^{d-1}}
\left\vert A_{j}\right\vert <q\left\vert c_{0}\right\vert +%
{\displaystyle\sum\limits_{j=1}^{d-1}}
\left\vert A_{j}\right\vert \leq%
{\displaystyle\sum\limits_{j=0}^{d-1}}
\left\vert A_{j}\right\vert .
\]
This also ends the proof of Theorem \ref{HRP} (ii), when $\alpha$ is an
algebraic integer, because by Theorem \ref{Dominant} (ii) we may choose the
polynomial $C$ so that $c_{d}=1,$ and the Euclidean division of any element
$Q$ $\in\mathbb{Z}[x]$ by $C$ gives that $Q(\alpha)=A_{0}+\cdot\cdot
\cdot+A_{d-1}\alpha^{d-1}$ for some $(A_{0},$ $...,$ $A_{d-1})\in
\mathbb{Z}^{d}.$

Now, we use a simple induction on $\deg(R)$ to complete the proof of Theorem
\ref{HRP}. By the above, we have $R(\alpha)\in F[\alpha],$ when $\deg(R)\leq
d-1.$ Let
\[
R(x)=A_{0}+A_{1}x+\cdot\cdot\cdot+A_{D}x^{D}\in\mathbb{Z}[x],
\]
where $D\geq d,$ and suppose that $P(\alpha)\in F[\alpha]$ for all
$P\in\mathbb{Z}[x],$ with $\deg(P)<D.$ Since $\deg(A_{0})=0\leq d-1,$ the
relation (3) implies that
\[
A_{0}=\varepsilon+\alpha(a_{0}+\cdot\cdot\cdot+a_{d-1}\alpha^{d-1}),
\]
for some $\varepsilon\in F$ and $a_{j}$ $\in\mathbb{Z.}$ Hence,
\[
R(\alpha)=\varepsilon+\alpha((a_{0}+A_{1})+\cdot\cdot\cdot+(a_{D-1}%
+A_{D})\alpha^{D-1}),
\]
where $a_{d}=...=a_{D-1}=0,$ and the induction hypothesis, applied to the
polynomial $(a_{0}+A_{1})+\cdot\cdot\cdot+(a_{D-1}+A_{D})x^{D-1}\in
\mathbb{Z}[x],$ leads to the desired result.\hfill$\Box$

\bigskip

\textbf{Proof of Theorem \ref{Circle}. \ }Let $\alpha$ be an algebraic number,
whose conjugates $\alpha^{(1)},...,\alpha^{(\deg(\alpha))}$ lie on the unit
circle. Since Theorem \ref{Circle} is true when $\alpha$ is a root of unity,
suppose that $\alpha$ is not an algebraic integer and the leading coefficient
$c$ of its minimal polynomial $M_{\alpha}$ satisfies $c\geq2.$

\textit{Case} $m(\alpha)=\deg(\alpha)/2.$ Set $m:=m(\alpha)$ and let
$\alpha^{(1)},\dots,\alpha^{(m)}$ be $m$ conjugates of $\alpha$ which are
multiplicatively independent. Without loss of generality, we may assume that
$\operatorname{Im}(\alpha^{(j)})>0$ for all $j\in\{1,...,m\}.$ Then, the map
$\Phi$ defined, from the field $\mathbb{Q}(\alpha)$ into the ring
$\mathbb{C}^{m},$ by the relation
\[
\Phi(\beta)=(\beta^{(1)},\dots,\beta^{(m)}),
\]
where $\beta^{(j)}$ is the image of $\beta$ by the conjugate map which sends
$\alpha$ to $\alpha^{(j)}$ $\forall j\in\{1,...,m\},$ is an embedding. We
shall show that there exist two positive real numbers $B=B(\alpha)$ and
$R=R(\alpha),$ such that for any $\beta_{0}\in\mathbb{Z}[\alpha]$ there are
$N=N(\alpha,\beta_{0})$ elements $s_{1},...,s_{N}$ of set $[0,B]\cap
\mathbb{N},$ and a number $\gamma\in\mathcal{O}:=\mathbb{Z}[\alpha
]\cap\mathbb{Z}[1/\alpha]$ satisfying
\begin{equation}
\beta_{0}=(%
{\displaystyle\sum\limits_{j=1}^{N}}
s_{j}\alpha^{j-1})+\gamma\alpha^{N}\text{ and }\Vert\Phi(\gamma)\Vert\leq R,
\end{equation}
where $\Vert.\Vert$ is the sup norm on the vector space $\mathbb{C}^{m}.$
Indeed, in this case, $\alpha$ satisfies the relation (1) with a finite subset
of $\mathbb{Z}\cap\lbrack-\max\{B,h\},\max\{B,h\}],$ where
\[
h:=\max\{h(\gamma)\text{ }|\text{ }\gamma\in E\},
\]
$h(\gamma)$ is the greatest modulus of the coefficients of a fixed
representation of $\gamma$ in $\mathbb{Z}[\alpha],$ and the set
\[
E:=\{\gamma\in\mathcal{O}\ |\text{ }\Vert\Phi(\gamma)\Vert\leq R\},
\]
is finite by Lemma \ref{Int}. If
\[
\beta=a_{0}+\cdot\cdot\cdot+a_{n}\alpha^{n}%
\]
for some $n\in\mathbb{N}$ and $\{a_{0},...,a_{n}\}\subset\mathbb{Z},$ then the
Euclidean division of $a_{0}$ by $c$ gives that there is $d\in\{0,1,...,c-1\}$
such that $\beta\equiv d\operatorname{mod}\alpha,$ i. e., $(\beta-d)/\alpha
\in\mathbb{Z}[\alpha].$ Moreover, since $M_{\alpha}(0)=c,$ the number $d$ is
unique. Hence, the map
\[
T:\beta\mapsto(\beta-d)/\alpha,
\]
is well defined from $\mathbb{Z}[\alpha]$ into itself. Now, fix $\beta_{0}%
\in\mathbb{Z}[\alpha],$ and set
\[
\beta_{k}:=\alpha\beta_{k+1}+d_{k+1},
\]
where $k\in\mathbb{N},$ $\beta_{k+1}=T(\beta_{k})$ and $d_{k+1}\in
\{0,1,...,c-1\}.$ Then
\[
\beta_{k+1}=\frac{\beta_{0}}{\alpha^{k+1}}-\frac{d_{1}}{\alpha^{k+1}}%
-\cdot\cdot\cdot-\frac{d_{k+1}}{\alpha^{1}}.
\]
With the notation of Lemma \ref{KroA}, set $R:=(43K(\alpha,m,2\pi/5)+10)c.$
By Lemma \ref{KroA}, there is $l\in\mathbb{N}\cap\lbrack0,K]$ such that
$\left\vert \arg(\beta_{0}^{(j)}/(\alpha^{(j)})^{l})\right\vert \leq2\pi/5$
for $j=1,\dots,m.$ Select $d_{l+1}^{\ast}$ such that $5Kc\leq d_{l+1}^{\ast
}<(5K+1)c,$ and $\beta_{l}\equiv d_{l+1}^{\ast}\operatorname{mod}\alpha$. Let
$\beta_{l+1}^{\ast}:=(\beta_{l}-d_{l+1}^{\ast})/\alpha$.
Putting $r:=d_{l+1}^{\ast}/5$ and $z:=\beta_{0}^{(j)}/(\alpha^{(j)})^{l}$ in
Lemma \ref{L2}, we obtain
\[
\left\vert \beta_{l+1}^{\ast(j)}\right\vert =\left\vert \frac{\beta_{0}^{(j)}%
}{(\alpha^{(j)})^{l}}-\sum_{j=1}^{l}\frac{d_{j}}{(\alpha^{(j)})^{l-j+1}%
}-d_{l+1}^{\ast}\right\vert <|\beta_{0}^{(j)}|\leq\Vert\Phi(\beta_{0})\Vert,
\]
when $(37K+8)c\leq|\beta_{0}^{(j)}|.$ On the other hand, if $|\beta_{0}%
^{(j)}|<(37K+8)c,$ then
\[
\left\vert \beta_{l+1}^{\ast(j)}\right\vert \leq(43K+9)c<R.
\]
This implies
\[
\Vert\Phi(\beta_{l+1}^{\ast})\Vert<\max\{R,\Vert\Phi(\beta_{0})\Vert\}
\]
and
\[
\beta_{0}=(%
{\displaystyle\sum\limits_{j=1}^{l}}
d_{j}\alpha^{j-1})+d_{l+1}^{\ast}\alpha^{l}+\beta_{l+1}^{\ast}\alpha^{l+1}.
\]
So we have
\[
\beta_{l+1}^{\ast}\in\beta_{0}/\alpha^{l+1}+\mathbb{Z}[1/\alpha]\subset
\alpha^{u}\mathbb{Z}[1/\alpha]
\]
with $u=\max\{0,n-l-1\}$. Iterating this procedure, we obtain a sequence
$(\beta_{l(j)+1}^{\ast})_{j=1,2,\dots}$ with $l=l(1)$ and $\beta
_{l(j)+1}^{\ast}\in\mathbb{Z}[1/\alpha]\cap\mathbb{Z}[\alpha]$ for
sufficiently large $j$. From Lemma \ref{Int}, $\Phi(\mathcal{O})$ have no
accumulation points in $\mathbb{C}^{m},$ we obtain that $\alpha$ can written
\[
\beta_{0}=(%
{\displaystyle\sum\limits_{j=1}^{N}}
s_{j}\alpha^{j-1})+\gamma\alpha^{N},
\]
where $N\in\mathbb{N}^{\ast},$ $s_{j}\in\lbrack0,B]\cap\mathbb{N},$
$B:=(5K+1)c$ and $\gamma\in E.$ Hence, (4) is true and this completes the
proof of the first implication in Theorem \ref{Circle}.

It follows immediately, from the case above, that $\alpha$ satisfies the
height reducing property, when $\deg(\alpha)=2,$ as $m(\alpha)=\deg(\alpha)/2$
(in this case the constant $K$ is much smaller and one can make explicit the
height given by the above proof).

\textit{Case} $m(\alpha)=\deg(\alpha)/2-1.$ The proof is almost the same but
we use Lemma \ref{MinusOne} instead of Lemma \ref{KroA}.

We are left to show the case $m(\alpha)=1.$ From Lemma \ref{Dep2}, any two
distinct conjugates $\alpha_{l}$ and $\alpha_{j},$ of $\alpha,$ in the upper
half plane, satisfy $\alpha_{l}^{b}\alpha_{j}^{b}=1$ or $\alpha_{l}%
^{b}\overline{\alpha_{j}}^{b}=1$ for some positive rational integer $b$. In
both cases, $\alpha^{b}$ has less number of conjugates than $\alpha$. We can
iterate this discussion until we find an integer, say again $b,$ such that the
only other conjugate of $\alpha^{b}$ is $\overline{\alpha}^{b}$. Then
$\alpha^{b}$ is quadratic and so by the case $m(\alpha^{b})=\deg(\alpha
^{b})/2,$ there is a finite subset $F$ of $\mathbb{Z}$ such that
$\mathbb{Z}[\alpha^{b}]=F[\alpha^{b}];$ thus $\mathbb{Z}[\alpha]=F[\alpha]$,
since any sum of the form $%
{\displaystyle\sum\limits_{j=0}^{s}}
c_{j}\alpha^{j},$ where $c_{j}\in\mathbb{Z},$ may be written
\[%
{\displaystyle\sum\limits_{j=0}^{s}}
c_{jb}\alpha^{jb}+\alpha%
{\displaystyle\sum\limits_{j=0}^{s}}
c_{1+jb}\alpha^{jb}+\cdot\cdot\cdot+\alpha^{b-1}%
{\displaystyle\sum\limits_{j=0}^{s}}
c_{b-1+jb}\alpha^{jb},
\]
with $c_{j}=0$ when $j\geq s+1.$ $\Box$ \bigskip

\textbf{Proof of Theorem \ref{Dominant}. }A direct application of Rouch\'{e}'s
theorem gives that a polynomial $P\in\mathbb{C}[x],$ with $k-$th strictly
dominant term, has exactly $k$ roots with modulus less than 1. The same
argument applied, in this case, to the polynomial $x^{\deg(P)}P(1/x)$ shows
that $P$ has $(\deg(P)-k)$ roots outside the closed unit disk (see also [9, p.
225]); thus P has no roots on the unit circle.

Now, suppose that $\alpha$ is a root of a non-zero (resp., of a monic) integer
polynomial, say again $P(x)=c_{0}+c_{1}x+\cdot\cdot\cdot+c_{\deg(P)}%
x^{\deg(P)},$ such that
\[
\left\vert c_{k}\right\vert \geq%
{\displaystyle\sum\limits_{j\neq k}}
\left\vert c_{j}\right\vert ,
\]
for some $k\in\{0,...,\deg(P)\}.$ Then, $\alpha$ is an algebraic number
(resp., an algebraic integer), and by the above we have that the direct
implication in Theorem \ref{Dominant} (ii) is true, since the conjugates of
$\alpha$ are among the roots of $P.$ To show the direct implication of Theorem
\ref{Dominant} (i), notice first, by Lemma \ref{Order}, that $\alpha$ is root
of unity, when it has a conjugate lying on the unit circle. Assume that
$\alpha$ is not a root of unity (so $\alpha$ has no conjugates on the unit
circle) and consider the polynomial
\[
P_{n}(x)=P(x)+(\varepsilon/n)x^{k},
\]
where $n\in\mathbb{N}^{\ast}$ and $\varepsilon=\mathrm{sgn}(c_{k}).$ Also, by
the above the polynomial $P_{n}$ has exactly $k$ roots inside the unit disk.
Let $\beta_{1,n},...,$ $\beta_{\deg(P),n}$ be the roots of $P_{n},$ and let
$\beta$ be a root of $P$ \ with modulus less than 1. Then, $\left\vert
P_{n}(\beta)\right\vert =\left\vert \beta^{k}/n\right\vert <1/n$ and so
$\lim_{n\rightarrow\infty}P_{n}(\beta)=0.$ It follows by the equation
$\lim_{n\rightarrow\infty}%
{\displaystyle\prod\limits_{1\leq j\leq\deg(P)}}
(\beta-\beta_{j,n})=0,$ that there is a subsequence of some sequence
$(\beta_{j_{0},n})_{n\geq1},$ where $j_{0}$ is fixed in $\{1,...,\deg(P)\},$
which converges to $\beta.$ Hence, $P$ has at most $k$ distinct roots with
modulus less than 1, and so $\alpha$ has at most $k$ conjugates inside the
unit disk, since its minimal polynomial is separable.

To prove the other implications in Theorem \ref{Dominant}, consider an
algebraic number (resp., an algebraic integer), say again $\alpha,$ having
$l\geq0$ conjugates with modulus less than 1 and no conjugates on the unit
circle. Then, by Lemma \ref{L4}, we see that there is $N\in\mathbb{N}^{\ast}$
such that the polynomial $Q(x):=%
{\displaystyle\prod\limits_{1\leq j\leq d}}
(x-\alpha_{j}^{N}),$ where $\alpha_{1},...,\alpha_{d}$ are the conjugates of
$\alpha,$ has an $l-$th strictly dominant term. Moreover, since $Q(x)\in
\mathbb{Q}[x],$ there is $v\in\mathbb{N}^{\ast}$ such that $vQ(x)\in
\mathbb{Z}[x],$ and so $\alpha$ is a root of the integer polynomial
$R(x)=vQ(x^{N})$ (resp., since $Q(x)\in\mathbb{Z}[x],$ $\alpha$ is a root of
the monic integer polynomial $R(x)=Q(x^{N}))$ with an $l-$th strictly dominant
term. Now, let $k\in\mathbb{N}\cap\lbrack l,\infty\lbrack.$ Then, $\alpha$ has
at most $k$ conjugates inside the unit disk, and is a root of the polynomial
\[%
{\displaystyle\sum\limits_{j=0}^{k-l-1}}
c_{j}^{\prime}x^{j}+x^{k-l}R(x),
\]
where $c_{j}^{\prime}=0$ for all $j\in\{0,...,k-l-1\},$ with $k-$th strictly
dominant term; this ends the proof of Theorem \ref{Dominant} (ii). Finally
notice when $\alpha$ is an $N-$th root of unity, then it is a root of the
monic integer polynomial $x^{2N+k}+(B-1)x^{N+k}-Bx^{k},$ where $B\in
\mathbb{N}^{\ast}$ and $k\in\mathbb{N,}$ with $k-$th dominant term, and this
completes the proof of Theorem \ref{Dominant} (i). \hfill$\Box$\bigskip
\ \bigskip

\begin{center}
\textbf{Appendix.}
\end{center}

Continuing Remark \ref{comp}, we describe briefly a practical method to study
multiplicative dependence of $\alpha_{i}$'s, by using Lemma 3.7 of \cite{deW}.
Put $\theta_{i}=\log\alpha_{i}$ for $i=1,\dots,m$ and $\theta_{m+1}=2\pi$.
Choose a large constant $C$. In this case, it seems enough to take $C=B^{m+2}$
where $B$ is the maximum of constants appearing in Lemma 4.1 of \cite{Wal}.
Apply LLL algorithm for the lattice generated by the following $m+1$ vectors:

\begin{eqnarray*}
&&(1,0,\dots,0,0,\lfloor C\theta_{1}\rfloor)\\
&&(0,1,0,\dots,0,\lfloor C\theta_{2}\rfloor)\\
&& \qquad \qquad \vdots\\
&&(0,0,0,\dots,1,\lfloor C\theta_{m}\rfloor)\\
&&(0,0,\dots,0,\lfloor C\theta_{m+1}\rfloor)\\
\end{eqnarray*}
where the notation $\lfloor.\rfloor$ designates the integer part function.

Using Proposition 1.11 of \cite{LLL}, if the first vector $v$ found by LLL
satisfies
\[
\Vert v\Vert>2^{m/2}\sqrt{(m^{2}+5m+4)}B,
\]
then $\alpha_{1},\dots,\alpha_{m}$ are multiplicatively independent since we
can choose $\delta$ in Lemma 3.7 of de Weger \cite{deW} as large as possible.
If this inequality does not hold, then the first vector $v=(k_{1}%
,\dots,k_{m+1})$ becomes small and it is highly possible that it gives a
multiplicative dependence $\prod_{j=1}^{m}\alpha_{j}^{k_{j}}=1$. We check the
validity by rigorous symbolic computation.

Hereafter we present some numerical results on the multiplicative dependency
of $\alpha$. It suggests that $m(\alpha)<\deg(\alpha)/2$ rarely happens.

Let us fix an even degree $d$ and a leading coefficient $c\geq2$. We are
interested in the number of primitive irreducible reciprocal polynomials of
degree $d,$ with the leading coefficient $c,$ whose all roots have modulus
one. Further if there is a positive rational integer $b$ such that
$\deg(\alpha^{b})<\deg(\alpha)$, then we can reduce the problem to lower
degree. By Lemma \ref{Dep2}, this occurs when and only when there are two
distinct multiplicatively dependent conjugates of $\alpha$ which are not
complex conjugates. We call this $\alpha$ \textit{power-reducible}. For e.g.,
$\alpha$ is power-reducible if the minimal polynomial $M_{\alpha}$ of $\alpha$
has a form $g(x^{m})$ for some rational integer $m\geq2$ and some polynomial
$g$. We wish to exclude power-reducible cases to obtain non trivial examples.
If $\deg(\alpha)\geq4$ and $m(\alpha)=1$ then $\alpha$ is certainly
power-reducible by Lemma \ref{Dep2}. The first non trivial case holds when
$d=6$ and $m(\alpha)=2$.

Put
\[
T_{n}^{\ast}(y)=%
\begin{cases}
2T_{n}(y/2) & n=1,2,\dots\\
1 & n=0
\end{cases}
\]
where $T_{n}(x)$ is the $n$-th Chebyshev polynomial of the 1-st kind. Fix a
positive rational integer $h$. To produce polynomials whose all roots are of
modulus one, we search integer polynomials
\[
g(y)=\sum_{j=0}^{d/2}c_{j}T_{d/2-j}^{\ast}(y)
\]
with $c_{0}=c$ and $\left\vert c_{j}\right\vert \leq h$ for all $j$. The
reciprocal polynomial
\[
c_{d/2}x^{d/2}+\sum_{j=0}^{d/2-1}c_{j}(x^{j}+x^{d-j})
\]
has $d$ roots on the unit circle if and only if $g(y)=0$ has $d/2$ real roots
in $[-2,2]$. We pick out such polynomials and check multiplicative dependence
by the method in Remark \ref{comp}. The result is shown in Table
\ref{Dependency} for $c=2$ and $c=3$.

We explain Table \ref{Dependency} by examples. 
Hereafter 
the index of complex roots in the upper half plane is sorted by real parts. 
For $(d,c,h)=(6,2,50)$, among
$1030301$ polynomials there are $287$ polynomials whose all roots are of modulus
one. Within them there are 62 primitive irreducible ones. There remain 58
polynomials which do not have the form $g(x^{m})$ with $m\geq2$. Finally using
the method of Remark \ref{comp}, we find 8 polynomials with $m(\alpha
)<\deg(\alpha)/2$. All of them satisfies $m(\alpha)=\deg(\alpha)/2-1$. For
e.g., $2-2x+3x^{2}-2x^{3}+3x^{4}-2x^{5}+2x^{6}$ gives $\alpha_{1}\alpha
_{2}^{-1}\alpha_{3}=\sqrt{-1}.$ 
For $(d,c,h)=(8,2,12)$, the above sieving
process does not suffice, because there are $16-10=6$ power-reducible
polynomials which does not have the form $g(x^{m})$ with $m\geq2$. For e.g,
let $\alpha$ be a root of
\[
2+4x+2x^{2}-4x^{3}-7x^{4}-4x^{5}+2x^{6}+4x^{7}+2x^{8}.
\]
Then $\alpha^{8}$ is a root of $16+8x+x^{2}+8x^{3}+16x^{4}$. The remaining 10
polynomials satisfy $m(\alpha)=\deg(\alpha)/2-1$.

We did not find any example which is not covered by Theorem \ref{Circle} for
degree not greater than $10$. Thus height reducing property is valid in this
search range of $c$ and $h$.

However in degrees $12$ and $16$, we find cases with
\[
m(\alpha)=\deg(\alpha)/2-2\text{ or }m(\alpha)=\deg(\alpha)/2-3.
\]
Such cases form pairs $\pm\alpha$ and we shall present one representative in
each pair.

\textit{Case $m(\alpha)=\deg(\alpha)/2-2$.}
\[
2 + 4 x + 4 x^{2} + 2 x^{3} + x^{4} + x^{8} + 2 x^{9} + 4 x^{10} + 4 x^{11} +
2 x^{12}
\]
whose dependencies are generated by $\alpha_{1}=\alpha_{4}\alpha_{5}^{-1}$ and
$\alpha_{2}=\alpha_{3}\alpha_{6}^{-1}$.
\[
3 - 3 x + x^{2} + x^{3} - 2 x^{4} + 2 x^{5} - x^{6} + 2 x^{7} - 2 x^{8} +
x^{9} + x^{10} - 3 x^{11} + 3 x^{12}%
\]
gives $\alpha_{1}\alpha_{6}/\alpha_{2}=\alpha_{3}\alpha_{5}/\alpha_{4}%
=\frac{1+\sqrt{-3}}2$.
\[
3 + 3 x^{2} - x^{4} - 2 x^{5} - 3 x^{6} - 2 x^{7} - x^{8} + 3 x^{10} + 3
x^{12}
\]
gives $\alpha_{1}\alpha_{3}/\alpha_{4}=\alpha_{2}\alpha_{5}/\alpha_{6}=-1$.
For degree $16$,
\[
2-2 x-x^{2}+x^{3}+x^{4}-2 x^{6}+x^{7}+x^{8}+x^{9}-2 x^{10}+x^{12}%
+x^{13}-x^{14}-2 x^{15}+2 x^{16}
\]
gives generating dependencies: $\alpha_{1}\alpha_{3}/(\alpha_{4}\alpha
_{8})=\alpha_{2}\alpha_{5}\alpha_{7}/\alpha_{6}=-1$. Adapting the idea of
Lemma \ref{MinusOne} simultaneously to two multiplicative dependences, we can
prove height reducing property for these 4 polynomials, by solving 9 systems
of inequalities.

\textit{Case $m(\alpha)=\deg(\alpha)/2-3$.}
\[
2+4x+4x^{2}+3x^{3}+3x^{4}+2x^{5}+x^{6}+2x^{7}+3x^{8}+3x^{9}+4x^{10}%
+4x^{11}+2x^{12}%
\]
gives $\alpha_{2}\alpha_{3}\alpha_{4}=\alpha_{1}\alpha_{3}\alpha_{5}=1$ and
$\alpha_{4}=\alpha_{5}\alpha_{6}$.
\[
3-3x^{2}+2x^{3}+3x^{4}-x^{6}+3x^{8}+2x^{9}-3x^{10}+3x^{12}%
\]
gives $\alpha_{3}\alpha_{4}/\alpha_{1}=\alpha_{3}\alpha_{5}/\alpha_{2}%
=\alpha_{2}\alpha_{6}/\alpha_{1}=1$. We are not able to show height reducing
property for these last two polynomials so far.

\begin{table}[h]%
\begin{tabular}
[c]{|r|r|r|r|r|r|r|r|r|r|r|r|r|}\hline
d & c & h & poly & circle & irred & prim & non $x^{m}$ & dep & npr & $-1$ & $-2$ &
-3\\\hline
6&	2&	50&	1030301& 287& 71&	62&	58&	8&	8&	8& 0 & 0\\
6&	3&	50&	1030301& 805&	325&	318&	310&	22&	22&	22& 0 & 0\\
8 & 2 & 12 & 390625 & 1069 & 210 & 200 & 182 & 16 & 10 & 10 & 0 & 0\\
8 & 3 & 12 & 390625 & 3991 & 1565 & 1558 & 1502 & 42 & 40 & 40 & 0 & 0\\
10 & 2 & 6 & 371293 & 2931 & 518 & 516 & 512 & 8 & 8 & 8 & 0 & 0\\
10 & 3 & 6 & 372193 & 13244 & 5640 & 5638 & 5630 & 72 & 72 & 72 & 0 & 0\\
12 & 2 & 4 & 531441 & 6557 & 1386 & 1380 & 1310 & 32 & 24 & 20 & 2 & 2\\
12 & 3 & 4 & 531441 & 33202 & 15858 & 15852 & 15620 & 98 & 90 & 84 & 4 & 2\\
14 & 2 & 3 & 823543 & 12185 & 2510 & 2510 & 2506 & 12 & 12 & 12 & 0 & 0\\
14 & 3 & 3 & 823543 & 70951 & 37548 & 37548 & 37544 & 120 & 120 & 120 & 0 &
0\\
16 & 2 & 2 & 390625 & 15143 & 3940 & 3934 & 3828 & 34 & 32 & 30 & 2 &
0\\\hline
\end{tabular}
\caption{Multiplicative Dependency}%
\label{Dependency}%
\end{table}

\begin{itemize}
\item $d$: degree of $\alpha$

\item $c$: the leading coefficient of the minimal polynomial $M_{\alpha}$.

\item $h$: the maximum modulus of the coefficients of $M_{\alpha}$.

\item $poly$: number of polynomials.

\item $circle$: number of polynomials whose all roots have modulus one.

\item $irred$: number of irreducible polynomials in {\it circle}.

\item $prim$: number of primitive polynomials in {\it irred}.

\item {\it non} $x^{m}$: number of polynomials satisfying $M_{\alpha}(x)\neq
g(x^{m})$ in {\it prim}.

\item $dep$: number of multiplicatively dependent cases among {\it non} $x^m$.

\item $npr$: number of non-power reducible polynomials in {\it dep}.

\item $-1$: number of polynomials with $m(\alpha)=\deg(\alpha)/2-1$ in {\it npr}.

\item $-2$: number of polynomials with $m(\alpha)=\deg(\alpha)/2-2$ in {\it npr}.

\item $-3$: number of polynomials with $m(\alpha)=\deg(\alpha)/2-3$ in {\it npr}.
\end{itemize}

\textbf{Acknowledgment.} We would like to express our gratitude to Professor
Attila Peth\H o for his advices concerning the practical computation of
$m(\alpha)$ in Remark \ref{comp} and Appendix.

\bigskip

Department of mathematics, Faculty of sciences, Niigata University, Ikarashi
2-8050, \ Niigata 950-2181, Japan

\bigskip

Department of mathematics and informatics, Larbi Ben M'hidi University, Oum El
Bouaghi \ \ 04000, \ Algeria

\bigskip

\textit{E-mail address}: akiyama@math.sc.niigata-u.ac.jp

\smallskip

\textit{E-mail address}: toufikzaimi@yahoo.com
\end{document}